\documentclass[reqno]{amsart}
\usepackage{amsmath,amssymb}
\usepackage{relsize}
\usepackage{graphicx,color}
\usepackage[all]{xy}
\newenvironment{psmallmatrix}
{\left(\begin{smallmatrix}}
	{\end{smallmatrix}\right)}

\theoremstyle{plain}
\newtheorem{introtheorem}{Theorem}
\newtheorem{introcorollary}[introtheorem]{Corollary}
\newtheorem{theorem}{Theorem}[section]

\theoremstyle{definition}
\newtheorem{definition}[theorem]{Definition}
\newtheorem{example}[theorem]{Example}

\theoremstyle{remark}
\newtheorem{remark}[theorem]{Remark}

\def\E{{\mathcal E}}

\def\D{{\mathcal D}}
\def\C{{\mathcal C}}

\def\P{{\mathcal P}}

\def\A{{\mathcal A}}

\def\cat0{\mathrm{cat}_0}

\def\dim{\mathrm{dim}}

\def\aut{\mathrm{aut}}

\begin{document}

\title[]
{Realising  a  finite group as a subgroup of a  product of two  groups of permutation matrices}

\author{Mahmoud Benkhalifa}
\address{Department of Mathematics. Faculty of  Sciences, University of Sharjah. Sharjah, United Arab Emirates}

\email{mbenkhalifa@sharjah.ac.ae}



\begin{abstract}  In this paper we prove that any finite group of order $n$  can be viewed as the group  of the solutions of a certain matrix equation $XB=BY$, where the unknowns $X,Y$ are two permutation matrices of order $n$ and $(1+k)n+2 $ respectively and where  $k\in \Bbb N$ is given by Cayley's theorem.  Moreover, we show that   $G$ is isomorphic to   a certain   subgroup   formed by  permutation matrices  of order $(1+k)n$ obtained  by permuting  all the rows of the identity matrix $I_{(1+k)n}$.
\end{abstract}
\maketitle
\section{Introduction}
Let $\P(n)$ denote   the group of permutation matrices of degree $n$. For   a given  matrix $B$,  let us consider the group  $\Omega_B$  of  the pairs $(X,Y)\in \P(n)\times \P(m)$ which are   solutions of the matrix equation $XB=BY$. Obviously,  $\Omega_B$ is finite group as  $\P(n)$ and $\P(m)$ are finite and it is worth noting that if $\lambda\in\Bbb Q$ and $(X,Y)\in \Omega_B$, then  the pair $(\lambda X,\lambda Y)$ needs not be in  $\Omega_B$  although that we have  $(\lambda X)B=B(\lambda Y)$  because  $\lambda X,\lambda Y$ are not permutation matrices for  $\lambda\neq 1$. 

A subgroup $H$ of $\P(n)$ is called \textit{realisable} if each element $M\in H$ is  obtained  by permuting  the rows of the identity matrix $I_n$ using a permutation  $\tau\in S_n$ satisfying $\tau(i)\neq i$ for all $1\leq i\leq n$. Here $S_n$ denotes  the symmetric group of order $n$.

Recall that, by   Cayley's theorem,   any finite group $G$ of order $n$ is isomorphic  to realisable  subgroup, denoted by $\C_G$, of  $P(n)$ via the map 
\begin{equation*}
G=\{g_1,\dots,g_n\}\rightarrow S_{n}\cong P(n)\,\,\,\,\,,\,\,\,\,\,\,g_{j}\mapsto \sigma _{j}=\left(
\begin{array}{cccc}
g_1 & g_2 & \ldots  & g_n \\
g_j & \sigma _{j}(g_2) & \ldots  & \sigma _{j}(g_n)%
\end{array}%
\right)\leftrightarrow M_{j}
\end{equation*}
where $M_j$  is the  matrix obtained  by permuting all  the rows of the identity matrix $I_n$ using $\sigma_j$.

Following the idea developed in \cite{Benk3} and inspired by the works done in \cite{B0,B10} regarding the so-called Kahn’s realisability problem of groups
( see \cite{Benk1, k} for more details), this paper is devoted to  answer the question whether a given  finite group  $G$  can occur as a group on the form  $\Omega_{B}$ and whether $G$ can be embedded in $\P(m)$, where $m>n$,  as a realisable subgroup. For this purpose we shall assign to $G$  a  matrix $B_G$  and a realisable  subgroup $\A_G$ of $\P((1+k)n+2) $,   where the number  $k$ is given by the decomposition of the permutation $\sigma_{2}$ into product of disjoint cycles, i.e,.  $\sigma_{2}=\tau_1\tau_2\dots\tau_k$ and we shall define   $\Omega_{B_G}$ as a certain   subgroup of $\A_G\times\C_G$.

\noindent The group $\A_G$ and the matrix $B_G$ are defined using  the framework of rational homotopy theory \cite{FHT}   and the ideas  developed  in \cite{Benk10,Benk3}.  More precisely,  $\A_G$ is defined in terms  of the cohomology of a certain a free  commutative cochain $\Bbb Q$-algebra  associated with the group $G$ and $B_G$ is related to its differential.

\smallskip

Thus,  in this paper we  establish the following result
 \begin{introtheorem}
	For any  finite group  $G$ of order $n$, there exists a  matrix $B_G$  such that $G$ is isomorphic to the group  $\Omega_{B_G}$ of the solutions of the matrix equation $XB_G=B_GY$, where the unknowns $X,Y$ are two permutation matrices belonging to the groups $\A_G$ and $\C_G$ respectively.
\end{introtheorem}
As a corollary   we derive
 \begin{introcorollary}
	Any  finite group  $G$ of order $n$ is isomorphic to  a realisable  subgroup   of  $\P((1+k)n)$.
\end{introcorollary}

\section{Main results}
\subsection{Definition of the group $\A_G$}
Let us start by recalling the main construction  in \cite{Benk3} on which this work is based. Indeed, let $G=\left\{ g_{1},g_{2},\ldots ,g_{n}\right\} $ be a finite  group of
order $n$ and let $S_n$ the symmetric group. By Cayley's theorem  there is a monomorphism 
\begin{equation*}\label{8}
\Psi:G\to S_{n}\,\,\,\,\,\,\,\,\,,\,\,\,\,\,\,\,\,\,\,\,\,\,\,\,\,g_{j}\mapsto \sigma _{j}:g_{k}\longrightarrow
g_{j}g_{k}\,\,\,\,\,\,\,\,,\,\,\,\,\,\,\,\,\,\,1\leq k\leq n
\end{equation*}

 For $2\leq j\leq n $, write
$\sigma _{j}=\left(
\begin{array}{cccc}
1 & 2 & \ldots  & n \\
j & \sigma _{j}(2) & \ldots  & \sigma _{j}(n)%
\end{array}%
\right)$ and let:
\begin{equation*}\label{55}
\sigma _{2}=\Big(1\,2\,\sigma _{2}(2)\ldots,\sigma^{\kappa_{2}}_{2}(2)\Big)\Big(i_{1}\sigma _{2}(i_{1})\,\ldots\,\sigma^{\kappa_{i_{1}}}_{2}(i_{1})\Big)\ldots\Big(i_{k}\,\sigma _{2}(i_{k})\ldots\sigma^{\kappa_{i_{k}}}_{2}(i_{k})\Big)
\end{equation*}
the decomposition  of $\sigma _{2}$ into a  product of cycles. 

\medskip
Recall that in \cite{Benk3} we constructed a free  commutative cochain $\Bbb Q$-algebra  $$\big(\Lambda(x_{1},x_{2},y_{1},y_{2},y_{3},\{z_{j},w_{j}\}_{g_i\in G}),\partial\big)$$ 
where the degrees of the elements in this graded algebra are  $$|x_{1}|=8\,\,\,\,\,\,\,\,,\,\,\,\,\,\,\,\,\,|x_{2}|=10\,\,\,\,\,\,\,\,,\,\,\,\,\,\,\,\,\, |w_{j}|=40$$ and where the differential is given by:
\begin{eqnarray}\label{28}
\partial(x_{1})&\hspace{-2mm}=&\hspace{-2mm}\partial(x_{2})=\partial(w_{j})=0\,,\,\,\,\,\,\partial(y_{1})=x^{3}_{1}x_{2},\,\,\,\,\,\partial(y_{2})=x^{2}_{1}x^{2}_{2},\,\,\,\,\,\partial(y_{3})=x_{1}x^{3}_{2}\nonumber\\
\partial (z_{j})&\hspace{-2mm}=&\hspace{-2mm}w_{j}^{3}+w_{j}w_{\sigma _{j+1}(1)}x^{4}_{2}+\overset{k}{{\underset{\tau=1}{\sum}}}w_{j}w_{\sigma_{j+1}(i_{\tau})}x^{4}_{2}+u+x^{15}_{1}\,\,\,\,,\,\,\,\,\,1\leq j\leq n-1\nonumber\\
\partial (z_{n})&\hspace{-2mm}=&\hspace{-2mm}w_{n}^{3}+w_{n}w_{1}x^{4}_{2}+\overset{k}{{\underset{\tau=1}{\sum}}}w_{n}w_{i_{\tau}}x^{4}_{2}+u+x^{15}_{1}
\end{eqnarray}
where $u=y_{1}y_{2}x^{4}_{1}x^{2}_{2}-y_{1}y_{3}x^{5}_{1}x_{2}+y_{2}y_{3}x^{6}_{1}$,   and we proved that
$$\E(\Lambda(x_{1},x_{2},y_{1},y_{2},y_{3},\{z_{j},w_{j}\}_{g_i\in G})\cong G$$ 
 where $\E(\Lambda(x_{1},x_{2},y_{1},y_{2},y_{3},\{z_{j},w_{j}\}_{1\leq j\leq n })$  denotes  the group  of  self homotopy  cochain equivalences   of  $\Lambda(x_{1},x_{2},y_{1},y_{2},y_{3},\{z_{j},w_{j}\}_{g_i\in G})$ (see \cite{Benk10,Benk3} for more details).
  
\bigskip 

Now let $V^{119}=\Bbb Q\{z_{1},\dots,z_{n}\}$ be the vector space spanned  by the  set $\{z_{1},\dots,z_{n}\}$.  Recall that $|z_{i}|=119$ for every $1\leq i\leq n$.  In  (\cite{Benk3},  Proposition 3.9),  it is shown that
$$\E(\Lambda(x_{1},x_{2},y_{1},y_{2},y_{3},\{z_{j},w_{j}\}_{g_i\in G })\cong\D^{119}_{40},$$
where $\D^{119}_{40}$ is the subgroup  of $\aut(V^{119})\times \E(\Lambda(x_{1},x_{2},y_{1},y_{2},y_{3},\{w_{j}\}_{g_i\in G })$ consisting of the couples $(\xi,[\alpha])$  making  the following
diagram commutes:
\begin{equation}\label{12}
\begin{picture}(300,90)(10,30)
\put(60,100){$V^{119}\hspace{1mm}\vector(1,0){186}\hspace{1mm}V^{119}$}
\put(68,76){\scriptsize $b$} \put(277,76){\scriptsize $b$}
\put(66,96){$\vector(0,-1){38}$} \put(275,96){$\vector(0,-1){38}$}
\put(175,103){\scriptsize $\xi$} \put(165,50){\scriptsize
	$H^{120}(\alpha)$} \put(60,46){$\Gamma_{G}^{120}
	\hspace{1mm}\vector(1,0){180}\hspace{1mm}\Gamma_{G}^{120}$}
\end{picture}
\end{equation}
where $\Gamma_{G}^{120}=H^{120}\big(\Lambda(x_{1},x_{2},y_{1},y_{2},y_{3},\{w_{j}\}_{g_i\in G}\big)$ and where  $b$ is defined by
\begin{equation}\label{3}
b(z_{i})=\widehat{\partial(z_{i})}\,\,\,\,\,\,\,,\,\,\,\,\,\,\,\,\,\,\,\,\,1\leq j\leq n
\end{equation}
Here  $\widehat{\partial(z_{i})}$ is  the cohomology class of  $\partial(z_{i})$ in $H^{120}\big(\Lambda(x_{1},x_{2},y_{1},y_{2},y_{3},\{w_{j}\}_{g_i\in G}\big).$

\noindent Moreover, it is shown that if  $(\xi,[\alpha])\in\D^{119}_{40}$, then there exists a unique permutation 
\begin{equation*}\label{5}
\sigma _{s}=\left(
\begin{array}{cccc}
1 & 2 & \ldots  & n \\
s & \sigma _{s}(2) & \ldots  & \sigma _{s}(n)%
\end{array}%
\right)
\end{equation*}
such that 
\begin{equation}\label{6}
\xi(z_{j})=z_{\sigma_{s}(j)}\,\,,\,\,\alpha(w_{j})=w_{\sigma_s(j)}\,\,,\,\,\alpha=id\,\,\,,\,\,\,\text{ on }\,\,x_{1},x_{2},y_{1},y_{2},y_{3}.
\end{equation}
 Thus, there is an isomorphism $\Psi:\D^{119}_{40}\to G$  defined by $\Psi((\xi,[\alpha]))=g_{s}$, 
 where the element $g_{s}$ corresponds to the permutation $\sigma _{s}$, given in (\ref{5}), via Cayley's theorem.
 
 \noindent Set $u=y_{1}y_{2}x^{4}_{1}x^{2}_{2}-y_{1}y_{3}x^{5}_{1}x_{2}+y_{2}y_{3}x^{6}_{1}$. As  the following set of generators 
\begin{equation}\label{20}
 \Sigma=\Big\{ w_{1}^{3}\,\,;\,\,\ldots\,\,;\,\,w_{n}^{3}\,\,;\,\,w_{j}w_{\sigma_{j+1}(1)}x^{4}_{2}\,\,;\,\,w_{j}w_{\sigma_{j+1}(i_{\tau})}x^{4}_{2}\,\,;\,\,u\,\,;\,\,x^{15}_{1}\Big\}
\end{equation}
 where  $1\leq j\leq n$ and $1\leq \tau\leq k$, is linearly independent in the vector space  $$\Gamma_{G}^{120}=H^{120}\big(\Lambda(x_{1},x_{2},y_{1},y_{2},y_{3},\{w_{j}\}_{i\in G}\big)$$ 
 it follows that   $\Sigma$  can be chosen, according the formulas (\ref{28}) and (\ref{3}), as   a basis  for the vector space  $b(V^{119})\subseteq\Gamma_{G}^{120}$. Notice that 
 \begin{equation}\label{9}
 \dim\,b(V^{119})=cardinal (\Sigma)=(1+k)n+2
 \end{equation}
 
 \noindent Thus,  if $B_{G}$  denotes the  the matrix of order $((1+k)n+2)\times n$  which is associated to the linear map $b$ defined in (\ref{3}) with respects to the basis is $\Sigma$, then we can write 
 $$
 B_{G}=
 \begin{bmatrix} 
 I_{n} \\
 M \\
 D 
 \end{bmatrix}
 \hspace{1cm}\text{ where }\hspace{1cm}
 D =\begin{bmatrix} 
 1 & 1 &\dots &1\\
 1 & 1 &\dots &1
 \end{bmatrix},
 $$
 where the matrix $M=\begin{bmatrix} 
 m_{ij} 
 \end{bmatrix}$ is defined by
 $$m_{ij} =
 \begin{cases}
 1,\,\,\,\,\text{ if } i\in\{\sigma_{j+1}(1),\sigma_{j+1}(i_{1}),\dots,\sigma_{j+1}(i_{k})\},\\
 0,\,\,\,\,\text{ otherwise}.
 \end{cases}
 $$

\noindent  Consequently, taking into construction (\ref{6}),   the matrices associated to the linear maps $\xi$ and the restriction of the linear map   $H^{120}(\alpha)$ to $b(V^{119})$, given in the diagram (\ref{12}) and corresponding to the element $(\xi,[\alpha])\in\D^{119}_{40}$, can be written, respectively,  as 
\begin{equation}\label{7}
C_{g_{s}}=\sigma_{s}{I_{n}}\,\,\,\,\,\,\,\,\,\,\,\,\,\,\,\,\,\,,\,\,\,\,\,\,\,\,\,\,\,\,\,\,\,\,\,\,\,\,\,A_{g_{s}}=
\begin{bmatrix} 
\sigma_{s}I_{n}& 0& 0\\
0&\widetilde{A}_{g_{s}}& 0 \\
0& 0&I_{2}
\end{bmatrix},
\end{equation}
where 
$$\sigma_{s}{I_{n}}=\begin{bmatrix} 
c_{i,j} 
\end{bmatrix}_{1\leq i,j\leq (k+1)n}\hspace{3mm}\text{ ,}\,\,\, c_{i,j} =
\begin{cases}
1,\,\,\,\,\text{ if } i
=\sigma_s(j)\\
0,\,\,\,\,\text{ otherwise}
\end{cases},$$
and where
$$\widetilde{A}_{g_{s}}=\begin{bmatrix} 
a_{n+i,n+j} 
\end{bmatrix}_{1\leq i,j\leq (k+1)n}\hspace{3mm}\text{ ,}\,\,\, a_{n+i,n+j} =
\begin{cases}
1,\,\,\,\,\text{ if } i
=\sigma_s(j)\\
0,\,\,\,\,\text{ otherwise}
\end{cases}.$$
Here $\sigma_{s}$ is the permutation   corresponds to  $g _{s}$ via Cayley's theorem.
\begin{remark}\label{r1}
From (\ref{7}), it is clear to see that $A_{g_{s}}$ is a permutation matrix. Recall that the commutativity of the diagram (\ref{12})  implies that 
\begin{equation}\label{2}
A_{g_s}B_G=B_GC_{g_s}\,\,\,\,\,\,\,\,\,\,\,,\,\,\,\,\,\,\,\,\,\,\,\,\forall g_s\in G
\end{equation}
\end{remark}
\begin{definition}\label{d1}
Let  $G=\left\{ g_{1},\ldots ,g_{n}\right\}$ be a group, we define  the following two sets
$$\A_G=\{A_{g_{s}}\,\,,\,\,g_{s}\in G\}\,\,\,\,\,\,\,\,\,\,,\,\,\,\,\,\,\,\,\,\,\,\,\Omega_G=\{(A_{g_{s}},C_{g_s})\in\A_G\times\C_G \,\,,\,\,g_{s}\in G\}$$ 	
\end{definition}	
\begin{theorem}\label{p2}
The sets   $\A_G$ and $\Omega_G$	are  groups isomorphic to $G$.
\end{theorem}	
\begin{proof}
	First let us prove that  $\A_G$  is a group. Indeed, let  $A_{g_s},A_{g_r}\in \A_G$. By (\ref{2}) there exist two  matrices $C_{g_s},C_{g_r}$ such that
	$$A_{g_{s}}B_{G}=B_{G}C_{g_{s}}\,\,\,\,\,\,\,\,\,\,\,\,\,\,\,\,\text{ and }\,\,\,\,\,\,\,\,\,\,\,\,\,\,\,\,A_{g_{r}}B_{G}=B_{G}C_{g_{r}}$$
	therefore $A_{g_{s}}A_{g_{r}}B_G=A_{g_{s}}B_GC_{g_{r}}=B_GC_{g_{s}}C_{g_{r}}$, it follows that $A_{g_{s}}A_{g_{r}}\in \A_G$. Here we use the that fact that   
\begin{equation}\label{11}
	A_{g_{s}}A_{g_{r}}=A_{g_{s}g_{r}}\,\,\,\,\,\,\,\,\,\,\,\,\,\,\,\,\text{ and }\,\,\,\,\,\,\,\,\,\,\,\,\,\,\,\,C_{g_{s}}C_{g_{r}}=C_{g_{s}g_{r}}
\end{equation}
	 Next let  $A_{g_s}\in \A_G$. As $A_{g_s}B_G=B_GC_{g_s}$ and $A_{g_s},C_{g_s}$ are invertible, we deduce that $B_GC_{g_s}^{-1}=(A_{g_s})^{-1}B_G$ implying that $(A_{g_s})^{-1}\in \A_G$. Notice that $A_{g_s}^{-1}=A_{g_s^{-1}}$.
	
		\noindent Next, using the same arguments, it is easy to  check that the set $\Omega_G$ is a group. Finally,  it is clear that the two maps $\chi:G \to \A_G$ and $\varphi:G \to \Omega_G$,  defined by $\chi(g_{s})=A_{g_{s}}$ and $\varphi(g_{s})=(A_{g_{s}},C_{g_s})$ respectively,  are  isomorphisms of groups.
		\end{proof}
\subsection{Realisable subgroups}
\begin{definition}\label{d2}
	A subgroup $H$ of $\P(n)$ is called \textit{realisable} if each element $M\in H$ is  obtained  by permuting  the rows of the identity matrix $I_n$ using a permutation  $\tau\in S_n$ satisfying $\tau(i)\neq i$ for all $1\leq i\leq n$. 
	\end{definition}
Let  $G=\left\{ g_{1},\ldots ,g_{n}\right\}$ be a group. Based on  the formula (\ref{7}) , let us define the following matrix 
\begin{equation}\label{10}
M_{g_{s}}=
\begin{bmatrix} 
\sigma_{s}I_{n}& 0\\
0&\widetilde{A}_{g_{s}} 
\end{bmatrix}\,\,\,\,\,\,\,\,\,\,\,\,\,\,,\,\,\,\,\,\,\,\,\,\,\,\,g_s\in G
\end{equation}
\begin{theorem}\label{tt1}
	If $H_G=\{M_{g_{s}}\,,\,\,g_s\in G\}$, then $H_G$ is a realisable subgroup of  $\P((1+k)n)$ isomorphic to $G$
\end{theorem}
\begin{proof}
According to 	the formula (\ref{7}), the matrix  $M_{g_{s}}$ is defined in terms of the permutation $\sigma_{s}$ corresponding to the element $g_s\in G$, so  it follows that $\sigma_{s}(i)\neq i$ for every $1\leq i\leq n$ implying that  $M_{g_{s}}\in\P((1+k)n)$. 

\noindent Taking into consideration the relation (\ref{11}), the map $G=\left\{ g_{1},\ldots ,g_{n}\right\}\to H_G$ which assign  $g_{s}\mapsto M_{g_{s}}$ is obviously an isomorphism of groups.
\end{proof}

\subsection{Examples}
In the following examples we illustrate our study by  determining all the groups introduced in this paper  for the cyclic group ${\Bbb Z_4}$ and the Klein group $\Bbb V$.
\begin{example}
	\label{e3}
	If $G=\Bbb Z_4$, then the monomorphism $\Bbb Z_{4}=\{ g_1,g_2,g_3,g_4\}\to S_4$ is given by 
	$$g_1\to id\,\,\,\,\,,\,\,\,\,\,g_2\leftrightarrow\sigma _{2}=(1234)\,\,\,\,\,,\,\,\,\,\,g_3\leftrightarrow\sigma_{3}=(13)(24)\,\,\,\,,\,\,\,\,\,g_4\leftrightarrow\sigma _{4}=(1432)$$	
	 therefore according to (\ref{28}) the model associated with   $\Bbb Z_4$ is $$\big(\Lambda(x_{1},x_{2},y_{1},y_{2},y_{3},w_{1},w_{2},w_{3},w_{4},z_{1},z_{2},z_{3},z_{4}),\partial\big)$$ 
	where $|x_{1}|=8,|x_{2}|=10, |w_{j}|=40$, and where the differential is given by
	\begin{eqnarray}
	\partial(x_{1})\hspace{-2mm}&=&\hspace{-2mm}\partial(x_{2})=\partial(w_{j})=0\,\,\,,\,\,\,\partial(y_{1})=x^{3}_{1}x_{2}\,\,,\,\,\,\partial(y_{2})=x^{2}_{1}x^{2}_{2}\,\,,\,\,\,\partial(y_{3})=x_{1}x^{3}_{2}\nonumber\\	
	\partial(z_{1})\hspace{-2mm}&=&\hspace{-2mm}w_{1}^{3}+w_{1}w_{2}x_{2}^4+y_{1}y_{2}x^{4}_{1}x^{2}_{2}-y_{1}y_{3}x^{5}_{1}x_{2}+y_{2}y_{3}x^{6}_{1}+x^{15}_{1}\nonumber\\
	\partial(z_{2})\hspace{-2mm}&=&\hspace{-2mm}w_{2}^{3}+w_{2}w_{3}x_{2}^4+y_{1}y_{2}x^{4}_{1}x^{2}_{2}-y_{1}y_{3}x^{5}_{1}x_{2}+y_{2}y_{3}x^{6}_{1}+x^{15}_{1}\nonumber\\
	\partial(z_{3})\hspace{-2mm}&=&\hspace{-2mm}w_{3}^{3}+w_{3}w_{4}x_{2}^4+y_{1}y_{2}x^{4}_{1}x^{2}_{2}-y_{1}y_{3}x^{5}_{1}x_{2}+y_{2}y_{3}x^{6}_{1}+x^{15}_{1}\nonumber\\
	\partial(z_{4})\hspace{-2mm}&=&\hspace{-2mm}w_{2}^{4}+w_{4}w_{1}x_{2}^4+y_{1}y_{2}x^{4}_{1}x^{2}_{2}-y_{1}y_{3}x^{5}_{1}x_{2}+y_{2}y_{3}x^{6}_{1}+x^{15}_{1}\nonumber
	\end{eqnarray} 
For the above construction it is clear that $V^{119}=\Bbb Q\{z_{1},z_{2},z_{3},z_{4}\}$ and by (\ref{20})	the base $\Sigma$ of the vector space $b(V^{119})$ is given 
	\begin{equation}\label{23}
	\Sigma=\Big\{ w_{1}^{3},w_{2}^{3},w_{3}^{3},w_{4}^{3},w_{1}w_{2}x^{4}_{2},w_{2}w_{3}x^{4}_{2},w_{3}w_{4}x^{4}_{2},w_{4}w_{1}x^{4}_{2},y_{1}y_{2}x^{4}_{1}x^{2}_{2}-y_{1}y_{3}x^{5}_{1}x_{2}+y_{2}y_{3}x^{6}_{1},x^{15}_{1}\Big\}
	\end{equation}
implying that the matrix $B_{\Bbb Z_4}$ associated with the linear map $b$, given in (\ref{3}), is 
$$B_{\Bbb Z_4}=
	\begin{psmallmatrix} 
	1 & 0  &0&0\\
	0 & 1  &0&0\\
	0 & 0 &1&0\\
	0 & 0  &0&1\\
	1 & 0  &0&0\\
    0& 1  &0&0\\
	0 & 0  &1&0\\
	0 & 0 &0&1\\
	1 & 1 &1&1\\
	1 & 1  &1&1
	\end{psmallmatrix}, 
$$
and we have   $\C_{\Bbb Z_4}=\Big\{I_4,C_{(1234)},C_{(13)(24)},C_{(1432)}\Big\}$ where
$$	C_{(1234)}=\begin{psmallmatrix} 
	0 & 0  &0 &1\\
	1 & 0  &0 &0\\
	0 & 1 &0 &0\\
	0 & 0 &1 &0
	\end{psmallmatrix},
	\hspace{.5cm}C_{(13)(24)}=\begin{psmallmatrix} 
	0 & 0 &1 &0\\
	0 & 0  &0 &1\\
	1 & 0 &0 &0\\
	0 & 1 &0&0
	\end{psmallmatrix}, 
	\hspace{.5cm}C_{(1432)}=\begin{psmallmatrix} 
	0 & 1&0&0\\
	0 & 0  &1 &0\\
	0 & 0 &0 &1\\
	1 & 0 &0&0
	\end{psmallmatrix}
	$$
For instance, $C_{(1234)}$ is simply the permutation matrix obtained by permuting the rows of $I_4$ using the permutation $(1234)$
	and likewise  $C_{(13)(24)}$ and $C_{(1432)}$

\noindent Next we have  $\A_{\Bbb Z_4}=\Big\{I_{10},A_{(1234)},A_{(13)(24)},A_{(1432)}\Big\}$ where
	$$A_{(1234)}=\begin{psmallmatrix} 
	0 & 0  &0&1 & 0 &0& 0  &0& 0   &0\\
	1 & 0 &0&0 & 0 &0 & 0  &0& 0  &0\\
	0 & 1 &0&0 & 0  &0& 0  &0& 0   &0 \\
	0 & 0 &1&0 & 0 &0& 0  &0& 0    &0 \\
	0 & 0 &0&0 & 0  &0& 0  &1& 0   &0 \\
	0 & 0 &0&0 & 1  &0& 0  &0& 0 &0\\
	0 & 0 &0&0 & 0  &1& 0  &0& 0  &0\\
	0 & 0 &0&0 & 0  &0& 1  &0& 0 &0 \\
	0 & 0 &0&0 & 0  &0& 0  &0& 1  &0 \\
	0 & 0 &0&0 & 0 &0& 0  &0& 0   &1
	\end{psmallmatrix}\,,\,\,\,
A_{(13)(24)}=\begin{psmallmatrix}
	0 & 0  &1&0 & 0 &0& 0  &0& 0   &0\\
	0& 0 &0&1 & 0 &0 & 0  &0& 0  &0\\
	1 & 0&0&0 & 0  &0& 0  &0& 0   &0 \\
	0 & 1 &0&0 & 0 &0& 0  &0& 0    &0 \\
	0 & 0 &0&0 & 0  &0& 1  &0& 0   &0 \\
	0 & 0 &0&0 & 0  &0& 0  &1& 0 &0\\
	0 & 0 &0&0 & 1  &0& 0  &0& 0  &0\\
	0 & 0 &0&0 & 0  &1&0 &0& 0 &0 \\
	0 & 0 &0&0 & 0  &0& 0  &0& 1  &0 \\
	0 & 0 &0&0 & 0 &0& 0  &0& 0   &1
	\end{psmallmatrix}\,,\,\,\,A_{(1432)}=\begin{psmallmatrix}
		0 & 1  &0&0 & 0 &0& 0  &0& 0   &0\\
		0& 0 &1&0 & 0 &0 & 0  &0& 0  &0\\
		0 & 0&0&1 & 0  &0& 0  &0& 0   &0 \\
		1 & 0&0&0 & 0 &0& 0  &0& 0    &0 \\
		0 & 0 &0&0 & 0  &1& 0  &0& 0   &0 \\
		0 & 0 &0&0 & 0  &0& 1  &0& 0 &0\\
		0 & 0 &0&0 & 0  &0& 0  &1& 0  &0\\
		0 & 0 &0&0 & 1  &0&0 &0& 0 &0 \\
		0 & 0 &0&0 & 0  &0& 0  &0& 1  &0 \\
		0 & 0 &0&0 & 0 &0& 0  &0& 0   &1
		\end{psmallmatrix}
		$$
		
		\noindent and  $H_{\Bbb Z_4}=\Big\{I_{8},M_{(1234)},M_{(13)(24)},M_{(1432)}\Big\}$, where
		$$M_{(1234)}=\begin{psmallmatrix} 
		0 & 0  &0&1 & 0 &0& 0  &0\\
		1 & 0 &0&0 & 0 &0 & 0  &0\\
		0 & 1 &0&0 & 0  &0& 0  &0 \\
		0 & 0 &1&0 & 0 &0& 0  &0 \\
		0 & 0 &0&0 & 0  &0& 0  &1 \\
		0 & 0 &0&0 & 1  &0& 0  &0\\
		0 & 0 &0&0 & 0  &1& 0  &0\\
		0 & 0 &0&0 & 0  &0& 1  &0
		\end{psmallmatrix}\,,\,\,\,
		M_{(13)(24)}=\begin{psmallmatrix}
		0 & 0  &1&0 & 0 &0& 0  &0\\
		0& 0 &0&1 & 0 &0 & 0  &0\\
		1 & 0&0&0 & 0  &0& 0  &0 \\
		0 & 1 &0&0 & 0 &0& 0  &0 \\
		0 & 0 &0&0 & 0  &0& 1  &0\\
		0 & 0 &0&0 & 0  &0& 0  &1\\
		0 & 0 &0&0 & 1  &0& 0  &0\\
		0 & 0 &0&0 & 0  &1&0 &0
		\end{psmallmatrix}\,,\,\,\,M_{(1432)}=\begin{psmallmatrix}
		0 & 1  &0&0 & 0 &0& 0  &0\\
		0& 0 &1&0 & 0 &0 & 0  &0\\
		0 & 0&0&1 & 0  &0& 0  &0 \\
		1 & 0&0&0 & 0 &0& 0  &0 \\
		0 & 0 &0&0 & 0  &1& 0  &0 \\
		0 & 0 &0&0 & 0  &0& 1  &0\\
		0 & 0 &0&0 & 0  &0& 0  &1\\
		0 & 0 &0&0 & 1  &0&0 &0
		\end{psmallmatrix}
		$$
Recall that $A_{(1234)}$ is the matrix associate to  the restriction of the linear map   $H^{120}(\alpha)$ to the vector space  $b(V^{119})$,  where the cochain map $\alpha$ is given (\ref{6}), with respects to the basis $\Sigma$ in (\ref{23}). Thus,  $A_{(1234)}$  is obtained by using the permutation $\sigma _{2}=(1234)$  as follows
$$w_{1}^{3}\mapsto w_{2}^{3}\,\,,\,\,w_{2}^{3}\mapsto w_{3}^{3}\,\,,\,\,w_{3}^{3}\mapsto w_{4}^{3}\,\,,\,\,w_{4}^{3}\mapsto w_{1}^{3}\,\,,\,\,w_{1}w_{2}x^{4}_{2}\mapsto w_{2}w_{3}x^{4}_{2}$$
$$w_{2}w_{3}x^{4}_{2}\mapsto w_{3}w_{4}x^{4}_{2}\,\,,\,\,w_{3}w_{4}x^{4}_{2}\mapsto w_{4}w_{1}x^{4}_{2}\,\,,\,\,w_{4}w_{1}x^{4}_{2}\mapsto w_{1}w_{4}x^{4}_{2}$$
$$y_{1}y_{2}x^{4}_{1}x^{2}_{2}-y_{1}y_{3}x^{5}_{1}x_{2}+y_{2}y_{3}x^{6}_{1}\mapsto y_{1}y_{2}x^{4}_{1}x^{2}_{2}-y_{1}y_{3}x^{5}_{1}x_{2}+y_{2}y_{3}x^{6}_{1}\,\,\,\,\,,\,\,\,\,x^{15}_{1}\mapsto x^{15}_{1}$$
and likewise we obtain the matrices $A_{(13)(24)}$ and $A_{(1432)}$.  Notice that matrices $M_{(1234)},M_{(13)(24)},M_{(1432)}$ are constructed from the matrices $A_{(1234)},A_{(13)(24)},A_{(1432)}$ using (\ref{11}) and   finally we have 
$$\Omega_{\Bbb Z_4}=\Big\{(I_{4},I_{12}),(A_{(1234)},C_{(1234)}),(A_{(13)(24)},C_{(13)(24)}),(A_{(1432)},C_{(1432)})\Big\}$$
It is also worth noting to point out that the group $M_{\Bbb Z_4}$, which is isomorphic to $\Bbb Z_4$  is  a realisable subgroup of the group of permutation matrices $\P(8)$.
\end{example}
\begin{example}
		\label{e2}
	In this example we  use the same analysis and computation as in the example  (\ref{e3}), but we omit  all the details, to  determine the two groups $\A_{\Bbb V},\C_{\Bbb V},H(\Bbb V)$ and $\Omega_{\Bbb V}$  for the Klein group $\Bbb V$. Indeed,   the monomorphism $\Bbb V=\{ g_1,g_2,g_3,g_4\}\leftrightarrow S_4$ is given by $g_2\leftrightarrow (12)(34)\,\,,\,g_3\leftrightarrow(13)(24)\,\,,\,g_4\to(14)(23)$, so  the model associated to  $\Bbb V$ is $\big(\Lambda(x_{1},x_{2},y_{1},y_{2},y_{3},w_{1},w_{2},w_{3},w_{4},z_{1},z_{2},z_{3},z_{4}),\partial\big)$
	where $|x_{1}|=8,|x_{2}|=10, |w_{j}|=40$, and where the differential is given by
	\begin{eqnarray}
\partial(x_{1})&=&\partial(x_{2})=\partial(w_{j})=0\,\,\,,\,\,\,\partial(y_{1})=x^{3}_{1}x_{2}\,\,,\,\,\,\partial(y_{2})=x^{2}_{1}x^{2}_{2}\,\,,\,\,\,\partial(y_{3})=x_{1}x^{3}_{2}\nonumber\\	
\partial(z_{1})&=&w_{1}^{3}+w_{1}w_{2}x_{2}^4+w_{1}w_{4}x_{2}^4+y_{1}y_{2}x^{4}_{1}x^{2}_{2}-y_{1}y_{3}x^{5}_{1}x_{2}+y_{2}y_{3}x^{6}_{1}+x^{15}_{1}\nonumber\\
\partial(z_{2})&=&w_{2}^{3}+w_{2}w_{3}x_{2}^4+w_{2}w_{1}x_{1}^4+y_{1}y_{2}x^{4}_{1}x^{2}_{2}-y_{1}y_{3}x^{5}_{1}x_{2}+y_{2}y_{3}x^{6}_{1}+x^{15}_{1}\nonumber\\
\partial(z_{3})&=&w_{3}^{3}+w_{3}w_{4}x_{2}^4+w_{3}w_{2}x_{2}^4+y_{1}y_{2}x^{4}_{1}x^{2}_{2}-y_{1}y_{3}x^{5}_{1}x_{2}+y_{2}y_{3}x^{6}_{1}+x^{15}_{1}\nonumber\\
\partial(z_{4})&=&w_{4}^{3}+w_{4}w_{1}x_{2}^4+w_{4}w_{3}x_{2}^4+y_{1}y_{2}x^{4}_{1}x^{2}_{2}-y_{1}y_{3}x^{5}_{1}x_{2}+y_{2}y_{3}x^{6}_{1}+x^{15}_{1}\nonumber
	\end{eqnarray}
If $u=y_{1}y_{2}x^{4}_{1}x^{2}_{2}-y_{1}y_{3}x^{5}_{1}x_{2}+y_{2}y_{3}x^{6}_{1}$, then the basis $\Sigma$ is given by
	\begin{equation}\label{023}
\Sigma=\Big\{ w_{1}^{3},w_{2}^{3},w_{3}^{3},w_{4}^{3},w_{1}w_{2}x^{4}_{2},w_{2}w_{3}x^{4}_{2},w_{3}w_{4}x^{4}_{2},w_{4}w_{1}x^{4}_{2},w_{1}w_{4}x^{4}_{2},w_{2}w_{1}x^{4}_{2},w_{3}w_{2}x^{4}_{2},w_{4}w_{3}x^{4}_{2},u,x^{15}_{1}\Big\}
\end{equation}	
	implying that $\dim\,b(V^{119})=14$ and the matrix $B_{\Bbb V}$ is 
	$$B_{\Bbb V}=\begin{psmallmatrix} 
	1 & 0  &0&0\\
	0 & 1  &0&0\\
	0 & 0 &1&0\\
	0 & 0  &0&1\\
	1 & 1  &0&0\\
	0& 0 &0&0\\
	1 & 0  &0&0\\
	0& 1  &0&0\\
	0 & 0  &1&1\\
	0 & 0 &1&0\\
	1 & 1 &1&1\\
	1 & 1  &1&1
	\end{psmallmatrix}
	$$
We have  $\C_{\Bbb V}=\Big\{I_4,C_{(12)(34)},C_{(13)(24)},C_{(14)(23)}\Big\}$, where
	$$
	C_{(12)(34)}=\begin{psmallmatrix} 
	0 & 1  &0 &0\\
	1 & 0  &0 &1\\
	0 & 0 &0 &0\\
	0 & 0 &1 &0
	\end{psmallmatrix}\,\,,
	\hspace{.2cm}C_{(13)(24)}=\begin{psmallmatrix} 
	0 & 0 &1 &0\\
	0 & 0  &0 &1\\
	1 & 0 &0 &0\\
	0 & 1 &0&0
	\end{psmallmatrix}\,\,,	\hspace{.2cm}C_{(14)(23)}=\begin{psmallmatrix} 
	0 & 0 &0&0\\
	0 & 0  &1 &0\\
	0 & 1 &0 &0\\
	1 & 0 &0&1
	\end{psmallmatrix}
	$$
	and   $\A_{\Bbb V}=\Big\{I_{14},A_{(12)(34)},A_{(13)(24)},A_{14)(23)}\Big\}$, where 
		$$	\hspace{-1.5mm} A_{(12)(34)}=\begin{psmallmatrix}
	0 & 1  &0&0 & 0  &0& 0  &0& 0  &0&0& 0  &0&0 \\
	1 & 0  &0&0 & 0  &0& 0  &0& 0  &0&0& 0  &0&0\\
	0 & 0  &0&1 & 0  &0& 0  &0& 0  &0&0& 0  &0&0 \\
	0 & 0  &1&0 & 0  &0& 0  &0& 0  &0&0& 0  &0&0 \\
	0 & 0  &0&0 & 0 &0& 0  &1& 0   &0&0& 0  &0&0 \\
	0 & 0 &0&0 & 0 &0 & 1   &0& 0  &0&0& 0  &0&0\\
	0 & 0 &0&0 & 0  &1& 0  &0& 0   &0&0& 0  &0&0 \\
	0 & 0 &0&0 & 1 &0& 0  &0& 0    &0&0& 0  &0&0 \\
	0 & 0 &0&0 & 0  &0& 0  &0& 0   &0&0& 1  &0&0 \\
	0 & 0 &0&0 & 0  &0& 0  &0& 0 &0&1& 0 &0&0 \\
	0 & 0 &0&0 & 0  &0& 0  &0& 0  &1&0& 0  &0&0 \\
	0 & 0 &0&0 & 0  &0& 0  &0& 1  &0&0& 0  &0&0 \\
	0 & 0 &0&0 & 0  &0& 0  &0& 0  &0&0& 0  &1&0 \\
	0 & 0 &0&0 & 0 &0& 0  &0& 0   &0&0& 0  &0&1
	\end{psmallmatrix}
	\,\,\,\,\,\,\,\,\,\,\,\,\,\,\,\,\,\,\,\,\,\,\,\,,\,\,\,\,\,\,\,\,\,\,\,\,\,\,\,\,
A_{(13)(24)}=\begin{psmallmatrix}
	0 & 0  &1&0 & 0  &0& 0  &0& 0  &0&0& 0  &0&0 \\
	0 & 0  &0&1 & 0  &0& 0  &0& 0  &0&0& 0  &0&0\\
	1 & 0  &0&0 & 0  &0& 0  &0& 0  &0&0& 0  &0&0 \\
	0 & 1  &0&0 & 0  &0& 0  &0& 0  &0&0& 0  &0&0 \\
	0 & 0  &0&0 &     0 &0& 0  &0      & 1   &0&0& 0  &0&0 \\
	0 & 0 &0&0 &      0 &0& 0  &0      & 0  &1&0& 0  &0&0\\
	0 & 0 &0&0 &      0 &0&0  &0      & 0   &0&1& 0  &0&0 \\
	0 & 0 &0&0 &      0 &0& 0  &0     & 0    &0&0& 1  &0&0 \\
	0 & 0 &0&0 &      1  &0& 0  &0&     0&0&0& 0           &0&0 \\
	0 & 0 &0&0 &      0  &1& 0  &0&     0 &0&0& 0          &0&0 \\
	0 & 0 &0&0 &      0  &0& 1  &0&     0&0&0& 0           &0&0 \\
	0 & 0 &0&0 &      0  &0& 0  &1&    0&0&0& 0         &0&0 \\
	0 & 0 &0&0 & 0  &0& 0  &0& 0  &0&0& 0  &1&0 \\
	0 & 0 &0&0 & 0 &0& 0  &0& 0   &0&0& 0  &0&1
	\end{psmallmatrix}
	$$
	$$A_{(14)(23)}=\begin{psmallmatrix}
	0 & 0  &0&1 & 0  &0& 0  &0& 0  &0&0& 0  &0&0 \\
	0 & 0  &1&0& 0  &0& 0  &0& 0  &0&0& 0  &0&0\\
	0 & 1 &0&0 & 0  &0& 0  &0& 0  &0&0& 0  &0&0 \\
	1 & 0 &0&0 & 0  &0& 0  &0& 0  &0&0& 0  &0&0 \\
	0 & 0  &0&0 &     0 &0& 0  &0      & 0   &0&0& 1  &0&0 \\
	0 & 0 &0&0 &      0 &0& 0  &0      & 0  &0&1& 0  &0&0\\
	0 & 0 &0&0 &      0 &0&0  &0      & 0   &1&0& 0  &0&0 \\
	0 & 0 &0&0 &      0 &0& 0  &0     & 1    &0&0& 0  &0&0 \\
	0 & 0 &0&0 &      0  &0& 0  &1&     0&0&0& 0           &0&0 \\
	0 & 0 &0&0 &      0  &0& 1  &0&     0 &0&0& 0          &0&0 \\
	0 & 0 &0&0 &      0  &1& 0  &0&     0&0&0& 0           &0&0 \\
	0 & 0 &0&0 &      1  &0& 0  &0&    0&0&0& 0        &0&0 \\
	0 & 0 &0&0 & 0  &0& 0  &0& 0  &0&0& 0  &1&0 \\
	0 & 0 &0&0 & 0 &0& 0  &0& 0   &0&0& 0  &0&1
	\end{psmallmatrix}
	$$
Next we have   $H_{\Bbb V}=\Big\{I_{12},M_{(12)(34)},M_{(13)(24)},M_{14)(23)}\Big\}$, where 
	$$	\hspace{-1.5mm} M_{(12)(34)}=\begin{psmallmatrix}
	0 & 1  &0&0 & 0  &0& 0  &0& 0  &0&0& 0  \\
	1 & 0  &0&0 & 0  &0& 0  &0& 0  &0&0& 0  \\
	0 & 0  &0&1 & 0  &0& 0  &0& 0  &0&0& 0   \\
	0 & 0  &1&0 & 0  &0& 0  &0& 0  &0&0& 0   \\
	0 & 0  &0&0 & 0 &0& 0  &1& 0   &0&0& 0   \\
	0 & 0 &0&0 & 0 &0 & 1   &0& 0  &0&0& 0 \\
	0 & 0 &0&0 & 0  &1& 0  &0& 0   &0&0& 0   \\
	0 & 0 &0&0 & 1 &0& 0  &0& 0    &0&0& 0   \\
	0 & 0 &0&0 & 0  &0& 0  &0& 0   &0&0& 1   \\
	0 & 0 &0&0 & 0  &0& 0  &0& 0 &0&1& 0 \\
	0 & 0 &0&0 & 0  &0& 0  &0& 0  &1&0& 0   \\
	0 & 0 &0&0 & 0  &0& 0  &0& 1  &0&0& 0  
	\end{psmallmatrix}
	\,\,\,\,\,\,\,\,\,\,\,\,\,\,\,\,\,\,\,\,\,\,\,\,,\,\,\,\,\,\,\,\,\,\,\,\,\,\,\,\,
	M_{(13)(24)}=\begin{psmallmatrix}
	0 & 0  &1&0 & 0  &0& 0  &0& 0  &0&0& 0   \\
	0 & 0  &0&1 & 0  &0& 0  &0& 0  &0&0& 0  \\
	1 & 0  &0&0 & 0  &0& 0  &0& 0  &0&0& 0   \\
	0 & 1  &0&0 & 0  &0& 0  &0& 0  &0&0& 0   \\
	0 & 0  &0&0 &     0 &0& 0  &0      & 1   &0&0& 0   \\
	0 & 0 &0&0 &      0 &0& 0  &0      & 0  &1&0& 0  \\
	0 & 0 &0&0 &      0 &0&0  &0      & 0   &0&1& 0   \\
	0 & 0 &0&0 &      0 &0& 0  &0     & 0    &0&0& 1   \\
	0 & 0 &0&0 &      1  &0& 0  &0&     0&0&0& 0            \\
	0 & 0 &0&0 &      0  &1& 0  &0&     0 &0&0& 0           \\
	0 & 0 &0&0 &      0  &0& 1  &0&     0&0&0& 0          \\
	0 & 0 &0&0 &      0  &0& 0  &1&    0&0&0& 0        
	\end{psmallmatrix}
	$$
	\vspace{0.5cm}
	$$M_{(14)(23)}=\begin{psmallmatrix}
	0 & 0  &0&1 & 0  &0& 0  &0& 0  &0&0& 0   \\
	0 & 0  &1&0& 0  &0& 0  &0& 0  &0&0& 0  \\
	0 & 1 &0&0 & 0  &0& 0  &0& 0  &0&0& 0  \\
	1 & 0 &0&0 & 0  &0& 0  &0& 0  &0&0& 0  \\
	0 & 0  &0&0 &     0 &0& 0  &0      & 0   &0&0& 1  \\
	0 & 0 &0&0 &      0 &0& 0  &0      & 0  &0&1& 0  \\
	0 & 0 &0&0 &      0 &0&0  &0      & 0   &1&0& 0   \\
	0 & 0 &0&0 &      0 &0& 0  &0     & 1    &0&0& 0  \\
	0 & 0 &0&0 &      0  &0& 0  &1&     0&0&0& 0           \\
	0 & 0 &0&0 &      0  &0& 1  &0&     0 &0&0& 0           \\
	0 & 0 &0&0 &      0  &1& 0  &0&     0&0&0& 0            \\
	0 & 0 &0&0 &      1  &0& 0  &0&    0&0&0& 0       
	\end{psmallmatrix}
	$$
Notice that  $H(\Bbb V)$ is a realizable  subgroup of $\P(12)$. 	Finally, we have $$\Omega_{\Bbb V}=\Big\{(I_{4},I_{14}),(A_{(12)(34)},C_{(12)(34)}),(A_{(13)(24)},C_{(13)(24)}),(A_{(14)(32)},C_{(14)(32)})\Big\}$$
\end{example}

\bibliographystyle{amsplain}

\end{document}